\documentclass[12pt,reqno]{amsart}
\usepackage{amsfonts}
\usepackage{amssymb,amsmath}
\usepackage{amscd}
\usepackage{hyperref}
\usepackage{graphicx}
\usepackage{epsfig}

\newtheorem{example}{Example}[section]

\numberwithin{equation}{section}
\begin{document}
\title[Pell Coding and Pell Decoding Methods]{Pell Coding and Pell Decoding
Methods with Some Applications}
\author[N. TA\c{S}]{N\.{I}HAL TA\c{S}*}
\address{Bal\i kesir University\\
Department of Mathematics\\
10145 Bal\i kesir, TURKEY}
\email{ nihaltas@balikesir.edu.tr}
\thanks{*Corresponding author: N. TA\c{S}\\
Bal\i kesir University, Department of Mathematics, 10145 Bal\i kesir, TURKEY%
\\
e-mail: nihaltas@balikesir.edu.tr}
\author[S. U\c{C}AR]{S\"{U}MEYRA U\c{C}AR}
\address{Bal\i kesir University\\
Department of Mathematics\\
10145 Bal\i kesir, TURKEY}
\email{sumeyraucar@balikesir.edu.tr}
\author[N. YILMAZ \"{O}ZG\"{U}R]{N\.{I}HAL YILMAZ \"{O}ZG\"{U}R}
\address{Bal\i kesir University\\
Department of Mathematics\\
10145 Bal\i kesir, TURKEY}
\email{nihal@balikesir.edu.tr}
\date{}

\begin{abstract}
We obtain a new coding and decoding method using the generalized Pell $(p,i)$%
-numbers. The relations among the code matrix elements, error detection and
correction have been established for this coding theory. We give two new
blocking algorithms using Pell numbers and generalized Pell $(p,i)$-numbers.
\end{abstract}

\maketitle

\section{Introduction}

\label{intro}

Recently, in \cite{Kilic-2009}, the generalized Pell $(p,i)$-numbers have
been defined by the following recurrence relations%
\begin{equation}
P_{p}^{(i)}(n)=2P_{p}^{(i)}(n-1)+P_{p}^{(i)}(n-p-1)\text{; }n>p+1,\text{ }%
0\leq i\leq p,\text{ }p=1,2,3,\ldots ,  \label{eqn1.1}
\end{equation}%
with the initial terms%
\begin{equation*}
P_{p}^{(i)}(1)=\cdots =P_{p}^{(i)}(i)=0\text{, }P_{p}^{(i)}(i+1)=\cdots
=P_{p}^{(i)}(p+1)=1\text{.}
\end{equation*}%
For $i=p=1$, the generalized Pell $(1,1)$-number corresponds to the $(n+1)$%
-th classical Pell number defined as
\begin{equation*}
P_{n+1}=2P_{n}+P_{n-1},\text{ }n\in
\mathbb{Z}
^{+}
\end{equation*}%
with the initial terms%
\begin{equation*}
P_{0}=0\text{, }P_{1}=1.
\end{equation*}%
For the basic properties of these numbers, one can see \cite{Kilic-2009} and
\cite{koshy-pell}.

It is known that $\gamma =1+\sqrt{2}$ and $\delta =1-\sqrt{2}$ are the roots
of the characteristic equation of the Pell recurrence relation $t^{2}-2t-1=0$%
. Using these roots it can be given the Binet formula for the Pell number by%
\begin{equation*}
P_{n}=\frac{\gamma ^{n}-\delta ^{n}}{\gamma -\delta }
\end{equation*}%
with%
\begin{equation*}
\underset{n\rightarrow \infty }{\lim }\frac{P_{n+1}}{P_{n}}=\gamma \text{.}
\end{equation*}

The Pell $P$-matrix of order $2$ is given by the following form:%
\begin{equation*}
P=\left[
\begin{array}{cc}
2 & 1 \\
1 & 0%
\end{array}%
\right] \text{.}
\end{equation*}%
The $n$-th power of the $P$-matrix and its determinant are given by%
\begin{equation}
P^{n}=\left[
\begin{array}{cc}
P_{n+1} & P_{n} \\
P_{n} & P_{n-1}%
\end{array}%
\right]   \label{pell power}
\end{equation}%
and%
\begin{equation*}
Det(P^{n})=P_{n+1}P_{n-1}-P_{n}^{2}=(-1)^{n}\text{.}
\end{equation*}

In \cite{Kilic-2009}, it was introduced the matrix $A$ in the following form%
\begin{equation}
A=\left[
\begin{array}{ccccc}
2 & 0 & \cdots & 0 & 1 \\
1 & 0 & \cdots & 0 & 0 \\
0 & 1 & \ddots & \vdots & 0 \\
\vdots & \cdots & \ddots & 0 & \vdots \\
0 & \cdots & 0 & 1 & 0%
\end{array}%
\right] _{(p+1)\times (p+1)}\text{.}  \label{eqn1.2}
\end{equation}%
and computed
\begin{equation}
{\footnotesize A^{n}=G_{n}=\left[
\begin{array}{ccccc}
P_{p}^{(p)}(n+p+1) & P_{p}^{(p)}(n+1) & P_{p}^{(p)}(n+2) & \cdots &
P_{p}^{(p)}(n+p) \\
P_{p}^{(p)}(n+p) & P_{p}^{(p)}(n) & P_{p}^{(p)}(n+1) & \cdots &
P_{p}^{(p)}(n+p-1) \\
\vdots & \vdots & \vdots & \vdots & \vdots \\
P_{p}^{(p)}(n+2) & P_{p}^{(p)}(n-p+2) & P_{p}^{(p)}(n-p+3) & \cdots &
P_{p}^{(p)}(n+1) \\
P_{p}^{(p)}(n+1) & P_{p}^{(p)}(n-p+1) & P_{p}^{(p)}(n-p+2) & \cdots &
P_{p}^{(p)}(n)%
\end{array}%
\right] _{(p+1)\times (p+1)}\text{.}}  \label{eqn1.3}
\end{equation}%
Using the matrices given in the equations (\ref{eqn1.2}) and (\ref{eqn1.3}),
we get%
\begin{equation}
Det(A^{n})=Det(G_{n})=(-1)^{n(p+2)}\text{,}  \label{eqn1.4}
\end{equation}%
for $n>0$.

Recently, Fibonacci coding theory has been introduced and studied in many
aspects (see \cite{basu}, \cite{prajat}, \cite{stakhov1999-2}, \cite{stakhov
2006}, \cite{Tarle} and \cite{Wang} for more details). For example, in \cite%
{stakhov 2006}, A. P. Stakhov presented a new coding theory using the
generalization of the Cassini formula for Fibonacci $p$-numbers and $Q_{p}$%
-matrices. Later B. Prasad developed a new coding and decoding method
followed from Lucas $p$ matrix \cite{prasad-lucas}. More recently it has
been obtained a new coding/decoding algorithm using Fibonacci $Q$-matrices
called as \textquotedblleft Fibonacci Blocking Algorithm\textquotedblright\
\cite{Tas}. In \cite{ucar}, using $R$-matrices and Lucas numbers, it has
been given \textquotedblleft Lucas Blocking Algorithm\textquotedblright\ and
the \textquotedblleft Minesweeper Model\textquotedblright\ related to
Fibonacci $Q^{n}$-matrices and $R$-matrices.

Motivated by the above studies, the main purpose of this paper is to develop
a new coding and decoding method using the generalized Pell $(p,i)$-numbers.
The relations among the code matrix elements, error detection and correction
have been established for this coding theory with $p=i=1$. As an
application, we give two new blocking algorithms using Pell and generalized
Pell $(p,i)$-numbers.

\section{Pell Coding and Decoding Method}

\label{sec:1} In this section, we present a new coding and decoding method
using the generalized Pell $(p,i)$-numbers. In our method, the nonsingular
square matrix $M$ with order $(p+1)$, where $p=1,2,3,\ldots $ corresponds to
our message matrix. We consider the matrix $G_{n}$ of order $(p+1)$ as
coding matrix and its inverse matrix $\left( G_{n}\right) ^{-1}$ as a
decoding matrix. We introduce Pell coding and Pell decoding with
transformations%
\begin{equation*}
M\times G_{n}=E
\end{equation*}%
and%
\begin{equation*}
E\times \left( G_{n}\right) ^{-1}=M\text{,}
\end{equation*}%
respectively, where $E$ is a code matrix.

Now we give an example of Pell coding and decoding method. Let $M$ be a
message matrix of the following form:%
\begin{equation*}
M=\left[
\begin{array}{cc}
m_{1} & m_{2} \\
m_{3} & m_{4}%
\end{array}%
\right] \text{,}
\end{equation*}%
where $m_{1},m_{2},m_{3},m_{4}$ are positive integers.

Let $p=i=1$ and $n=3$ in order to construct the coding matrix $G_{n}$:%
\begin{equation*}
G_{3}=\left[
\begin{array}{cc}
P_{1}^{(1)}(n+2) & P_{1}^{(1)}(n+1) \\
P_{1}^{(1)}(n+1) & P_{1}^{(1)}(n)%
\end{array}%
\right] =\left[
\begin{array}{cc}
P_{1}^{(1)}(5) & P_{1}^{(1)}(4) \\
P_{1}^{(1)}(4) & P_{1}^{(1)}(3)%
\end{array}%
\right] =\left[
\begin{array}{cc}
12 & 5 \\
5 & 2%
\end{array}%
\right] \text{.}
\end{equation*}%
Then we find inverse of $G_{3}$:%
\begin{equation*}
\left( G_{3}\right) ^{-1}=\frac{1}{12.2-5.5}\left[
\begin{array}{cc}
2 & -5 \\
-5 & 12%
\end{array}%
\right] =\left[
\begin{array}{cc}
-2 & 5 \\
5 & -12%
\end{array}%
\right] \text{.}
\end{equation*}%
Now, we calculate the code matrix $E$.%
\begin{equation*}
E=M\times G_{3}=\left[
\begin{array}{cc}
m_{1} & m_{2} \\
m_{3} & m_{4}%
\end{array}%
\right] \left[
\begin{array}{cc}
12 & 5 \\
5 & 2%
\end{array}%
\right] =\left[
\begin{array}{cc}
12m_{1}+5m_{2} & 5m_{1}+2m_{2} \\
12m_{3}+5m_{4} & 5m_{3}+2m_{4}%
\end{array}%
\right] \text{,}
\end{equation*}%
where $e_{1}=12m_{1}+5m_{2}$, $e_{2}=5m_{1}+2m_{2}$, $e_{3}=12m_{3}+5m_{4}$
and $e_{4}=5m_{3}+2m_{4}$.

Finally, the code matrix $E$ is sent to a channel, the message matrix $M$
can be obtained by decoding as the following way%
\begin{eqnarray*}
M &=&E\times \left( G_{3}\right) ^{-1}=\left[
\begin{array}{cc}
e_{1} & e_{2} \\
e_{3} & e_{4}%
\end{array}%
\right] \left[
\begin{array}{cc}
-2 & 5 \\
5 & -12%
\end{array}%
\right] \\
&=&\left[
\begin{array}{cc}
-2e_{1}+5e_{2} & 5e_{1}-12e_{2} \\
-2e_{3}+5e_{4} & 5e_{3}-12e_{4}%
\end{array}%
\right] =\left[
\begin{array}{cc}
m_{1} & m_{2} \\
m_{3} & m_{4}%
\end{array}%
\right] \text{.}
\end{eqnarray*}

Notice that the relation between the code matrix $E$ and the message matrix $%
M$ is%
\begin{equation*}
Det(E)=Det(M\times G_{n})=Det(M)\times Det(G_{n})=Det(M)\times (-1)^{n(p+2)}%
\text{.}
\end{equation*}

\section{The Relationships between the Code Matrix Elements for $p=i=1$}

\label{sec:2} We write $E$ and $M$ as follows:%
\begin{equation*}
E=M\times G_{n}=\left[
\begin{array}{cc}
m_{1} & m_{2} \\
m_{3} & m_{4}%
\end{array}%
\right] \left[
\begin{array}{cc}
P_{1}^{(1)}(n+2) & P_{1}^{(1)}(n+1) \\
P_{1}^{(1)}(n+1) & P_{1}^{(1)}(n)%
\end{array}%
\right] =\left[
\begin{array}{cc}
e_{1} & e_{2} \\
e_{3} & e_{4}%
\end{array}%
\right]
\end{equation*}%
and%
\begin{eqnarray*}
M &=&E\times \left( G_{n}\right) ^{-1}=\left[
\begin{array}{cc}
e_{1} & e_{2} \\
e_{3} & e_{4}%
\end{array}%
\right] \frac{1}{(-1)^{n(p+2)}}\left[
\begin{array}{cc}
P_{1}^{(1)}(n) & -P_{1}^{(1)}(n+1) \\
-P_{1}^{(1)}(n+1) & P_{1}^{(1)}(n+2)%
\end{array}%
\right]  \\
&=&\frac{1}{(-1)^{n(p+2)}}\left[
\begin{array}{cc}
e_{1}P_{1}^{(1)}(n)-e_{2}P_{1}^{(1)}(n+1) &
-e_{1}P_{1}^{(1)}(n+1)+e_{2}P_{1}^{(1)}(n+2) \\
e_{3}P_{1}^{(1)}(n)-e_{4}P_{1}^{(1)}(n+1) &
-e_{3}P_{1}^{(1)}(n+1)+e_{4}P_{1}^{(1)}(n+2)%
\end{array}%
\right] \text{,}
\end{eqnarray*}%
for $n=2k+1$. Because $m_{1},m_{2},m_{3},m_{4}$ are positive integers, we get%
\begin{equation}
m_{1}=-e_{1}P_{1}^{(1)}(n)+e_{2}P_{1}^{(1)}(n+1)>0\text{,}  \label{eqn3.1}
\end{equation}%
\begin{equation}
m_{2}=e_{1}P_{1}^{(1)}(n+1)-e_{2}P_{1}^{(1)}(n+2)>0\text{,}  \label{eqn3.2}
\end{equation}%
\begin{equation}
m_{3}=-e_{3}P_{1}^{(1)}(n)+e_{4}P_{1}^{(1)}(n+1)>0\text{,}  \label{eqn3.3}
\end{equation}%
\begin{equation}
m_{4}=e_{3}P_{1}^{(1)}(n+1)-e_{4}P_{1}^{(1)}(n+2)>0\text{.}  \label{eqn3.4}
\end{equation}%
Using (\ref{eqn3.1}) and (\ref{eqn3.2}), we find%
\begin{equation}
\frac{P_{1}^{(1)}(n+2)}{P_{1}^{(1)}(n+1)}<\frac{e_{1}}{e_{2}}<\frac{%
P_{1}^{(1)}(n+1)}{P_{1}^{(1)}(n)}\text{.}  \label{eqn3.5}
\end{equation}%
Using (\ref{eqn3.3}) and (\ref{eqn3.4}), we get%
\begin{equation}
\frac{P_{1}^{(1)}(n+2)}{P_{1}^{(1)}(n+1)}<\frac{e_{3}}{e_{4}}<\frac{%
P_{1}^{(1)}(n+1)}{P_{1}^{(1)}(n)}\text{.}  \label{eqn3.6}
\end{equation}%
From the inequalities (\ref{eqn3.5}) and (\ref{eqn3.6}), we obtain%
\begin{equation}
\frac{e_{1}}{e_{2}}\approx \gamma \text{, }\frac{e_{3}}{e_{4}}\approx \gamma
\text{.}  \label{eqn3.7}
\end{equation}%
For $n=2k$, we obtain the similar relations given in (\ref{eqn3.7}).

\section{Error Detection/Correction for Pell Coding and Decoding Method}

\label{sec:4}

Because of the reasons arising in the channel, some errors may be occur in
the code matrix $E.$ So we try to determine and correct these errors using
the properties of determinant in this process. Let $p=i=1$ and the message
matrix $M$ is given by
\begin{equation*}
M=\left[
\begin{array}{cc}
m_{1} & m_{2} \\
m_{3} & m_{4}%
\end{array}%
\right] .
\end{equation*}

We know that $Det(M)=m_{1}m_{4}-m_{2}m_{3},$ the code matrix $E=M\times
G_{n} $ and $Det(E)=Det(M)\times \left( -1\right) ^{n\left( p+2\right) }.$
From the relationship between determinants of $E$ and $M,$ if $n$ is an odd
number, we have%
\begin{equation*}
Det(E)=-Det(M)
\end{equation*}%
and if $n$ is an even number, we have
\begin{equation*}
Det(E)=Det(M).
\end{equation*}%
The new method of the error detection is an application of the $G_{n}$
matrix. The basic idea of this method depends on calculating the
determinants of $M$ and $E$. Comparing the determinants obtained from the
channel, the receiver can decide whether the code message $E$ is true or not.

Actually, we cannot determine which element of the code message is damaged.
In order to find damaged element, we suppose different cases such as single
error, two errors etc. Now, we consider the first case with a single error
in the code matrix $E$. We can easily obtain that there are four places
where single error appear in $E:$
\begin{eqnarray*}
&&(1)\left[
\begin{array}{cc}
t_{1} & e_{2} \\
e_{3} & e_{4}%
\end{array}%
\right] \text{,} \\
&&(2)\left[
\begin{array}{cc}
e_{1} & t_{2} \\
e_{3} & e_{4}%
\end{array}%
\right] \text{,} \\
&&(3)\left[
\begin{array}{cc}
e_{1} & e_{2} \\
t_{3} & e_{4}%
\end{array}%
\right] \text{, } \\
&&(4)\left[
\begin{array}{cc}
e_{1} & e_{2} \\
e_{3} & t_{4}%
\end{array}%
\right] \text{,}
\end{eqnarray*}%
where $t_{i}$ $\left( i\in \{1,2,3,4\}\right) $ is the damaged element. To
check the above four different cases, we can use the following relations$:$%
\begin{equation}
t_{1}e_{4}-e_{2}e_{3}=\left( -1\right) ^{n\left( p+2\right) }Det(M)\text{,}
\label{eqn41}
\end{equation}%
\begin{equation}
e_{1}e_{4}-t_{2}e_{3}=\left( -1\right) ^{n\left( p+2\right) }Det(M)\text{,}
\label{eqn42}
\end{equation}%
\begin{equation}
e_{1}e_{4}-t_{3}e_{2}=\left( -1\right) ^{n\left( p+2\right) }Det(M)\text{,}
\label{eqn43}
\end{equation}%
\begin{equation}
e_{1}t_{4}-e_{2}e_{3}=\left( -1\right) ^{n\left( p+2\right) }Det(M)\text{,}
\label{eqn44}
\end{equation}%
where a possible single error is the element $e_{1}$ $\left( \text{resp. }%
e_{2},e_{3}\text{ and }e_{4}\right) $ given in the relation $($\ref{eqn41}$)$
$\left( \text{resp. }(\text{\ref{eqn42}}),(\text{\ref{eqn43}})\text{ and }(%
\text{\ref{eqn44}})\right) .$

Using the above relations, we have
\begin{equation}
t_{1}=\frac{\left( -1\right) ^{n\left( p+2\right) }Det(M)+e_{2}e_{3}}{e_{4}}%
\text{,}  \label{eqn45}
\end{equation}%
\begin{equation}
t_{2}=\frac{\left( -1\right) ^{n\left( p+1\right) }Det(M)+e_{1}e_{4}}{e_{3}}%
\text{,}  \label{eqn46}
\end{equation}%
\begin{equation}
t_{3}=\frac{\left( -1\right) ^{n\left( p+1\right) }Det(M)+e_{1}e_{4}}{e_{2}}%
\text{,}  \label{eqn47}
\end{equation}%
\begin{equation}
t_{4}=\frac{\left( -1\right) ^{n\left( p+2\right) }Det(M)+e_{2}e_{3}}{e_{1}}%
\text{.}  \label{eqn48}
\end{equation}%
Since the elements of message matrix $M$ are positive integers, we should
find integer solutions of the equations from $($\ref{eqn45}$)$ to $($\ref%
{eqn48}$)$. If there are not integer solutions of these equations, we find
that our cases related to a single error is incorrect or an error can be
occurred in the checking element \textquotedblleft $Det(M)$%
\textquotedblright . If $Det(M)$ is incorrect, we use the relations given in
$($\ref{eqn3.7}$)$ to check a correctness of the code matrix $E$.

Similarly, we can check the cases with double errors in the code matrix $E$.
Let us consider the following case with double error in $E:$

\begin{equation}
\left[
\begin{array}{cc}
t_{1} & t_{2} \\
e_{3} & e_{4}%
\end{array}%
\right] ,  \label{eqn49}
\end{equation}%
where $t_{1},t_{2}$ are the damaged elements of $E.$ Using the relation
\begin{equation*}
Det(E)=\left( -1\right) ^{n\left( p+2\right) }Det(M),
\end{equation*}
we can write following equation for the matrix given in $($\ref{eqn49}$)$%
\begin{equation}
t_{1}e_{4}-t_{2}e_{3}=\left( -1\right) ^{n\left( p+2\right) }Det(M).
\label{eqn410}
\end{equation}%
Also, we know the following relation between $t_{1}$ and $t_{2}$%
\begin{equation}
t_{1}\approx \gamma t_{2}.  \label{eqn411}
\end{equation}%
It is clear that the equation $($\ref{eqn410}$)$ is a Diophantine equation.
Because there are many solutions of Diophantine equations, we should choose
the solutions $t_{1},t_{2}$ satisfying the above checking relation $($\ref%
{eqn411}$).$ Using the similar approach, we can correct the triple errors in
the code matrix $E$ such that
\begin{equation*}
\left[
\begin{array}{cc}
t_{1} & t_{2} \\
t_{3} & e_{4}%
\end{array}%
\right] ,
\end{equation*}%
where $t_{1},t_{2}$ and $t_{3}$ are damaged elements of $E$.

Consequently, our method depends on confirmation of various cases about
damaged elements of $E$ using the checking element $Det(M)$ and checking
relation $t_{1}\approx \gamma t_{2}.$ Our correctness method permits us to
correct $14$ cases among the $15$ cases, because our method is inadequate
for the case with four errors. So we can say that correction ability of our
method is $\frac{14}{15}=\%0,9333.$

\section{Blocking Methods As an Application of Pell Numbers and Generalized
Pell $(p,i)$-Numbers}

\label{sec:5}

In this section we introduce new coding/decoding algorithms using Pell and
generalized Pell $(p,i)$-numbers. We put our message in a matrix of even
size adding zero between two words and end of the message until we obtain
the size of the message matrix is even. Dividing the message square matrix $M
$ of size $2m$ into the block matrices, named $B_{i}$ ($1\leq i\leq m^{2}$),
of size $2\times 2$, from left to right, we can construct a new coding
method.

Now we explain the symbols of our coding method. Assume that matrices $B_{i}$
and $E_{i}$ are of the following forms:%
\begin{equation*}
B_{i}=\left[
\begin{array}{cc}
b_{1}^{i} & b_{2}^{i} \\
b_{3}^{i} & b_{4}^{i}%
\end{array}%
\right] \text{ and }E_{i}=\left[
\begin{array}{cc}
e_{1}^{i} & e_{2}^{i} \\
e_{3}^{i} & e_{4}^{i}%
\end{array}%
\right] \text{.}
\end{equation*}%
We use the matrix $P^{n}$ given in (\ref{pell power}) and rewrite the
elements of this matrix as $P^{n}=\left[
\begin{array}{cc}
p_{1} & p_{2} \\
p_{3} & p_{4}%
\end{array}%
\right] $. The number of the block matrices $B_{i}$ is denoted by $b$.
According to $b$, we choose the number $n$ as follows:%
\begin{equation*}
n=\left\{
\begin{array}{ccc}
3 & \text{,} & b\leq 3 \\
\left\lfloor \frac{b}{2}\right\rfloor & \text{,} & b>3%
\end{array}%
\right. \text{.}
\end{equation*}%
Using the chosen $n$, we write the following character table according to $%
mod29$ (this table can be extended according to the used characters in the
message matrix). We begin the \textquotedblleft $n$\textquotedblright\ for
the last character.

\begin{equation*}
\begin{tabular}{|c|c|c|c|c|c|c|c|c|c|}
\hline
A & B & C & D & E & F & G & H & I & J \\ \hline
$n+28$ & $n+27$ & $n+26$ & $n+25$ & $n+24$ & $n+23$ & $n+22$ & $n+21$ & $%
n+20 $ & $n+19$ \\ \hline
K & L & M & N & O & P & Q & R & S & T \\ \hline
$n+18$ & $n+17$ & $n+16$ & $n+15$ & $n+14$ & $n+13$ & $n+12$ & $n+11$ & $%
n+10 $ & $n+9$ \\ \hline
U & V & W & X & Y & Z & 0 & $:$ & $)$ &  \\ \hline
$n+8$ & $n+7$ & $n+6$ & $n+5$ & $n+4$ & $n+3$ & $n+2$ & $n+1$ & $n$ &  \\
\hline
\end{tabular}%
\end{equation*}

Now we explain the following new coding and decoding algorithms.

\textbf{Pell Blocking Algorithm}

\textbf{Coding Algorithm}

\textbf{Step 1.} Divide the matrix $M$ into blocks $B_{i}$ $\left( 1\leq
i\leq m^{2}\right) $.

\textbf{Step 2.} Choose $n$.

\textbf{Step 3. }Determine $b_{j}^{i}$ $\left( 1\leq j\leq 4\right) $.

\textbf{Step 4.} Compute $\det (B_{i})\rightarrow d_{i}$.

\textbf{Step 5.} Construct $K=\left[ d_{i},b_{k}^{i}\right] _{k\in
\{1,3,4\}} $.

\textbf{Step 6. }End of algorithm.

\textbf{Decoding Algorithm }

\textbf{Step 1.} Compute $P^{n}$.

\textbf{Step 2.} Determine $p_{j}$ $(1\leq j\leq 4)$.

\textbf{Step 3. }Compute $p_{1}b_{3}^{i}+p_{3}b_{4}^{i}\rightarrow e_{3}^{i}$
$\left( 1\leq i\leq m^{2}\right) $.

\textbf{Step 4.} Compute $p_{2}b_{3}^{i}+p_{4}b_{4}^{i}\rightarrow e_{4}^{i}$%
.

\textbf{Step 5.} Solve $(-1)^{n}\times
d_{i}=e_{4}^{i}(p_{1}b_{1}^{i}+p_{3}x_{i})-e_{3}^{i}(p_{2}b_{1}^{i}+p_{4}x_{i})
$.

\textbf{Step 6. }Substitute for $x_{i}=b_{2}^{i}$.

\textbf{Step 7.} Construct $B_{i}$.

\textbf{Step 8.} Construct $M$.

\textbf{Step 9.} End of algorithm.

In the following example we give an application of the above algorithm for $%
b>3$.

\begin{example}
\label{exm2} Let us consider the message matrix for the following message
text$:$%
\begin{equation*}
\text{\textquotedblleft MATH\ IS\ SWEET:)\textquotedblright }
\end{equation*}%
Using the message text, we get the following message matrix $M:$%
\begin{equation*}
M=\left[
\begin{array}{cccc}
M & A & T & H \\
0 & I & S & 0 \\
S & W & E & E \\
T & : & ) & 0%
\end{array}%
\right] _{4\times 4}.
\end{equation*}%
\textbf{Coding Algorithm:}

\textbf{Step 1. }We can divide the message text $M$ of size $4\times 4$ into
the matrices, named $B_{i}$ $\left( 1\leq i\leq 4\right) $, from left to
right, each of size is $2\times 2:$%
\begin{equation*}
B_{1}=\left[
\begin{array}{cc}
M & A \\
0 & I%
\end{array}%
\right] \text{, }B_{2}=\left[
\begin{array}{cc}
T & H \\
S & 0%
\end{array}%
\right] \text{, }B_{3}=\left[
\begin{array}{cc}
S & W \\
T & :%
\end{array}%
\right] \text{ and }B_{4}=\left[
\begin{array}{cc}
E & E \\
) & 0%
\end{array}%
\right] \text{.}
\end{equation*}

\textbf{Step 2.} Since $b=4>3$, we calculate $n=\left\lfloor \frac{b}{2}%
\right\rfloor =2$. For $n=2$, we use the following \textquotedblleft letter
table\textquotedblright\ for the message matrix $M:$%
\begin{equation*}
\begin{tabular}{|l|l|l|l|l|l|l|l|}
\hline
$M$ & $A$ & $T$ & $H$ & $0$ & $I$ & $S$ & $0$ \\ \hline
$18$ & $1$ & $11$ & $23$ & $4$ & $22$ & $12$ & $4$ \\ \hline
$S$ & $W$ & $E$ & $E$ & $T$ & $:$ & $)$ & $0$ \\ \hline
$12$ & $8$ & $26$ & $26$ & $11$ & $3$ & $2$ & $4$ \\ \hline
\end{tabular}%
\text{.}
\end{equation*}

\textbf{Step 3.} We have the elements of the blocks $B_{i}$ $\left( 1\leq
i\leq 4\right) $ as follows:%
\begin{equation*}
\begin{tabular}{|l|l|l|l|}
\hline
$b_{1}^{1}=18$ & $b_{2}^{1}=1$ & $b_{3}^{1}=4$ & $b_{4}^{1}=22$ \\ \hline
$b_{1}^{2}=11$ & $b_{2}^{2}=23$ & $b_{3}^{2}=12$ & $b_{4}^{2}=4$ \\ \hline
$b_{1}^{3}=12$ & $b_{2}^{3}=8$ & $b_{3}^{3}=11$ & $b_{4}^{3}=3$ \\ \hline
$b_{1}^{4}=26$ & $b_{2}^{4}=26$ & $b_{3}^{4}=2$ & $b_{4}^{4}=4$ \\ \hline
\end{tabular}%
.
\end{equation*}

\textbf{Step 4.} Now we calculate the determinants $d_{i}$ of the blocks $%
B_{i}:$%
\begin{equation*}
\begin{tabular}{|l|}
\hline
$d_{1}=\det (B_{1})=392$ \\ \hline
$d_{2}=\det (B_{2})=-232$ \\ \hline
$d_{3}=\det (B_{3})=-52$ \\ \hline
$d_{4}=\det (B_{4})=52$ \\ \hline
\end{tabular}%
.
\end{equation*}

\textbf{Step 5.} Using Step 3 and Step 4 we obtain the following matrix $K:$%
\begin{equation*}
K=\left[
\begin{array}{cccc}
392 & 18 & 4 & 22 \\
-232 & 11 & 12 & 4 \\
-48 & 12 & 11 & 3 \\
52 & 26 & 2 & 4%
\end{array}%
\right] .
\end{equation*}

\textbf{Step 6.} End of algorithm.\newline
\textbf{Decoding algorithm:}

\textbf{Step 1.} By $($\ref{pell power}$)$ we know that%
\begin{equation*}
P^{2}=\left[
\begin{array}{cc}
P_{3} & P_{2} \\
P_{2} & P_{1}%
\end{array}%
\right] =\left[
\begin{array}{cc}
5 & 2 \\
2 & 1%
\end{array}%
\right] \text{.}
\end{equation*}

\textbf{Step 2. }The elements of $P^{2}$ are denoted by%
\begin{equation*}
p_{1}=5\text{, }p_{2}=2\text{, }p_{3}=2\text{ and }p_{4}=1\text{.}
\end{equation*}

\textbf{Step 3.} We compute the elements $e_{3}^{i}$ to construct the matrix
$E_{i}:$%
\begin{equation*}
e_{3}^{1}=64\text{, }e_{3}^{2}=68\text{, }e_{3}^{3}=61\text{ and }%
e_{3}^{4}=18\text{.}
\end{equation*}

\textbf{Step 4. }We compute the elements $e_{4}^{i}$ to construct the matrix
$E_{i}:$%
\begin{equation*}
e_{4}^{1}=30\text{, }e_{4}^{2}=28\text{, }e_{4}^{3}=25\text{ and }e_{4}^{4}=8%
\text{.}
\end{equation*}

\textbf{Step 5.} We calculate the elements $x_{i}:$%
\begin{eqnarray*}
(-1)^{2}(392) &=&30(90+2x_{1})-64(36+x_{1}) \\
&\Rightarrow &x_{1}=1\text{.}
\end{eqnarray*}%
\begin{eqnarray*}
(-1)^{2}(-232) &=&28(55+2x_{2})-68(22+x_{2}) \\
&\Rightarrow &x_{2}=23\text{.}
\end{eqnarray*}%
\begin{eqnarray*}
(-1)^{2}(-52) &=&25(60+2x_{3})-61(24+x_{3}) \\
&\Rightarrow &x_{3}=8\text{.}
\end{eqnarray*}%
\begin{eqnarray*}
(-1)^{2}(52) &=&8(130+2x_{4})-18(52+x_{4}) \\
&\Rightarrow &x_{4}=26\text{.}
\end{eqnarray*}

\textbf{Step 6.} We rename $x_{i}$ as follows$:$%
\begin{equation*}
x_{1}=b_{2}^{1}=1\text{, }x_{2}=b_{2}^{2}=23\text{, }x_{3}=b_{2}^{3}=8\text{
and }x_{4}=b_{2}^{4}=26\text{.}
\end{equation*}

\textbf{Step 7. }We construct the block matrices $B_{i}:$%
\begin{equation*}
B_{1}=\left[
\begin{array}{cc}
18 & 1 \\
4 & 22%
\end{array}%
\right] \text{, }B_{2}=\left[
\begin{array}{cc}
11 & 13 \\
12 & 4%
\end{array}%
\right] \text{, }B_{3}=\left[
\begin{array}{cc}
12 & 8 \\
11 & 3%
\end{array}%
\right] \text{ and }B_{4}=\left[
\begin{array}{cc}
26 & 26 \\
2 & 4%
\end{array}%
\right] \text{.}
\end{equation*}

\textbf{Step 8.} We obtain the message matrix $M:$%
\begin{equation*}
M=\left[
\begin{array}{cccc}
18 & 1 & 11 & 23 \\
4 & 22 & 12 & 4 \\
12 & 8 & 26 & 26 \\
11 & 3 & 2 & 4%
\end{array}%
\right] =\left[
\begin{array}{cccc}
M & A & T & H \\
0 & I & S & 0 \\
S & W & E & E \\
T & : & ) & 0%
\end{array}%
\right] .
\end{equation*}

\textbf{Step 9.} End of algorithm.
\end{example}

Now we give an application of generalized Pell $(p,i)$-numbers using
blocking method. Assume that matrices $B_{i}$ and $E_{i}$ are of the
following forms:%
\begin{equation*}
B_{i}=\left[
\begin{array}{cc}
b_{1}^{i} & b_{2}^{i} \\
b_{3}^{i} & b_{4}^{i}%
\end{array}%
\right] \text{ and }E_{i}=\left[
\begin{array}{cc}
e_{1}^{i} & e_{2}^{i} \\
e_{3}^{i} & e_{4}^{i}%
\end{array}%
\right] \text{.}
\end{equation*}%
We use the matrix $G_{n}$ defined in (\ref{eqn1.3}) and rewrite the elements
of this matrix as $G_{n}=\left[
\begin{array}{cc}
g_{1} & g_{2} \\
g_{3} & g_{4}%
\end{array}%
\right] $. The number of the block matrices $B_{i}$ is denoted by $b$.
According to $b$, we choose the number $n$ as follows:
\begin{equation*}
n=p+2.
\end{equation*}

\textbf{Generalized Pell Blocking Algorithm}

\textbf{Coding Algorithm}

We consider coding algorithm like as the Pell Blocking coding algorithm.

\textbf{Decoding Algorithm }

\textbf{Step 1.} Compute $G_{n}$.

\textbf{Step 2.} Determine $g_{j}$ $(1\leq j\leq 4)$.

\textbf{Step 3. }Compute $g_{1}b_{3}^{i}+g_{3}b_{4}^{i}\rightarrow e_{3}^{i}$
$\left( 1\leq i\leq m^{2}\right) $.

\textbf{Step 4.} Compute $g_{2}b_{3}^{i}+g_{4}b_{4}^{i}\rightarrow e_{4}^{i}$%
.

\textbf{Step 5.} Solve $(-1)^{n(p+2)}\times
d_{i}=e_{4}^{i}(g_{1}b_{1}^{i}+g_{3}x_{i})-e_{3}^{i}(g_{2}b_{1}^{i}+g_{4}x_{i})
$.

\textbf{Step 6. }Substitute for $x_{i}=b_{2}^{i}$.

\textbf{Step 7.} Construct $B_{i}$.

\textbf{Step 8.} Construct $M$.

\textbf{Step 9.} End of algorithm.

We choose $p=i=1$ in the following example.

\begin{example}
\label{exm1} Let us consider the message matrix for the following message
text$:$%
\begin{equation*}
\text{\textquotedblleft HAPPY BIRTHDAY TO YOU:)\textquotedblright }
\end{equation*}%
Using the message text, we get the following message matrix%
\begin{equation*}
M=\left[
\begin{array}{cccccc}
H & A & P & P & Y & 0 \\
B & I & R & T & H & D \\
A & Y & 0 & T & O & 0 \\
Y & O & U & : & ) & 0 \\
0 & 0 & 0 & 0 & 0 & 0 \\
0 & 0 & 0 & 0 & 0 & 0%
\end{array}%
\right] _{6\times 6}.
\end{equation*}%
\textbf{Coding Algorithm:}

\textbf{Step 1. }We can divide the message text $M$ of size $6\times 6$ into
the matrices, named $B_{i}$ $\left( 1\leq i\leq 9\right) $, from left to
right, each of size is $2\times 2:$%
\begin{eqnarray*}
B_{1} &=&\left[
\begin{array}{cc}
H & A \\
B & I%
\end{array}%
\right] \text{, }B_{2}=\left[
\begin{array}{cc}
P & P \\
R & T%
\end{array}%
\right] \text{, }B_{3}=\left[
\begin{array}{cc}
Y & 0 \\
H & D%
\end{array}%
\right] \text{,} \\
B_{4} &=&\left[
\begin{array}{cc}
A & Y \\
Y & O%
\end{array}%
\right] \text{, }B_{5}=\left[
\begin{array}{cc}
0 & T \\
U & :%
\end{array}%
\right] \text{, }B_{6}=\left[
\begin{array}{cc}
O & 0 \\
) & 0%
\end{array}%
\right] \text{,} \\
B_{7} &=&\left[
\begin{array}{cc}
0 & 0 \\
0 & 0%
\end{array}%
\right] \text{, }B_{8}=\left[
\begin{array}{cc}
0 & 0 \\
0 & 0%
\end{array}%
\right] \text{, }B_{9}=\left[
\begin{array}{cc}
0 & 0 \\
0 & 0%
\end{array}%
\right] \text{.}
\end{eqnarray*}

\textbf{Step 2.} Since $p=1$, we calculate $n=3$. For $n=3$, we use the
following \textquotedblleft letter table\textquotedblright\ for the message
matrix $M:$%
\begin{equation*}
\begin{tabular}{|l|l|l|l|l|l|l|l|l|l|l|l|}
\hline
$H$ & $A$ & $P$ & $P$ & $Y$ & $0$ & $B$ & $I$ & $R$ & $T$ & $H$ & $D$ \\
\hline
$24$ & $2$ & $16$ & $16$ & $7$ & $5$ & $1$ & $23$ & $14$ & $12$ & $24$ & $28$
\\ \hline
$A$ & $Y$ & $0$ & $T$ & $O$ & $0$ & $Y$ & $O$ & $U$ & $:$ & $)$ & $0$ \\
\hline
$2$ & $7$ & $5$ & $12$ & $17$ & $5$ & $7$ & $17$ & $11$ & $4$ & $3$ & $5$ \\
\hline
$0$ & $0$ & $0$ & $0$ & $0$ & $0$ & $0$ & $0$ & $0$ & $0$ & $0$ & $0$ \\
\hline
$5$ & $5$ & $5$ & $5$ & $5$ & $5$ & $5$ & $5$ & $5$ & $5$ & $5$ & $5$ \\
\hline
\end{tabular}%
\text{.}
\end{equation*}

\textbf{Step 3.} We have the elements of the blocks $B_{i}$ $\left( 1\leq
i\leq 9\right) $ as follows:%
\begin{equation*}
\begin{tabular}{|l|l|l|l|}
\hline
$b_{1}^{1}=24$ & $b_{2}^{1}=2$ & $b_{3}^{1}=1$ & $b_{4}^{1}=23$ \\ \hline
$b_{1}^{2}=16$ & $b_{2}^{2}=16$ & $b_{3}^{2}=14$ & $b_{4}^{2}=12$ \\ \hline
$b_{1}^{3}=7$ & $b_{2}^{3}=5$ & $b_{3}^{3}=24$ & $b_{4}^{3}=28$ \\ \hline
$b_{1}^{4}=2$ & $b_{2}^{4}=7$ & $b_{3}^{4}=7$ & $b_{4}^{4}=17$ \\ \hline
$b_{1}^{5}=5$ & $b_{2}^{5}=12$ & $b_{3}^{5}=11$ & $b_{4}^{5}=4$ \\ \hline
$b_{1}^{6}=17$ & $b_{2}^{6}=5$ & $b_{3}^{6}=3$ & $b_{4}^{6}=5$ \\ \hline
$b_{1}^{7}=5$ & $b_{2}^{7}=5$ & $b_{3}^{7}=5$ & $b_{4}^{7}=5$ \\ \hline
$b_{1}^{8}=5$ & $b_{2}^{8}=5$ & $b_{3}^{8}=5$ & $b_{4}^{8}=5$ \\ \hline
$b_{1}^{9}=5$ & $b_{2}^{9}=5$ & $b_{3}^{9}=5$ & $b_{4}^{9}=5$ \\ \hline
\end{tabular}%
.
\end{equation*}

\textbf{Step 4.} Now we calculate the determinants $d_{i}$ of the blocks $%
B_{i}:$%
\begin{equation*}
\begin{tabular}{|l|}
\hline
$d_{1}=\det (B_{1})=550$ \\ \hline
$d_{2}=\det (B_{2})=-32$ \\ \hline
$d_{3}=\det (B_{3})=76$ \\ \hline
$d_{4}=\det (B_{4})=-15$ \\ \hline
$d_{5}=\det (B_{5})=-112$ \\ \hline
$d_{6}=\det (B_{6})=70$ \\ \hline
$d_{7}=\det (B_{7})=0$ \\ \hline
$d_{8}=\det (B_{8})=0$ \\ \hline
$d_{9}=\det (B_{9})=0$ \\ \hline
\end{tabular}%
.
\end{equation*}

\textbf{Step 5.} Using Step 3 and Step 4 we obtain the following matrix $K:$%
\begin{equation*}
K=\left[
\begin{array}{cccc}
550 & 24 & 1 & 23 \\
-32 & 16 & 14 & 12 \\
76 & 7 & 24 & 28 \\
-15 & 2 & 7 & 17 \\
-112 & 5 & 11 & 4 \\
70 & 17 & 3 & 5 \\
0 & 5 & 5 & 5 \\
0 & 5 & 5 & 5 \\
0 & 5 & 5 & 5%
\end{array}%
\right] .
\end{equation*}

\textbf{Step 6.} End of algorithm.\newline
\textbf{Decoding algorithm:}

\textbf{Step 1.} Using the equation $($\ref{eqn1.3}$)$ we have%
\begin{equation*}
G_{n}=\left[
\begin{array}{cc}
P_{p}^{(p)}\left( n+p+1\right) & P_{p}^{(p)}\left( n+1\right) \\
P_{p}^{(p)}\left( n+p\right) & P_{p}^{(p)}\left( n\right)%
\end{array}%
\right]
\end{equation*}%
and%
\begin{equation*}
G_{3}=\left[
\begin{array}{cc}
P_{1}^{(1)}\left( 5\right) & P_{1}^{(1)}\left( 4\right) \\
P_{1}^{(1)}\left( 4\right) & P_{1}^{(1)}\left( 3\right)%
\end{array}%
\right] =\left[
\begin{array}{cc}
12 & 5 \\
5 & 2%
\end{array}%
\right] .
\end{equation*}

\textbf{Step 2. }The elements of $G_{3}$ are denoted by
\begin{equation*}
g_{1}=12\text{, }g_{2}=5\text{, }g_{3}=5\text{ and }g_{4}=2.
\end{equation*}

\textbf{Step 3.} We compute the elements $e_{3}^{i}$ to construct the matrix
$E_{i}:$%
\begin{equation*}
e_{3}^{1}=127\text{, }e_{3}^{2}=228\text{, }e_{3}^{3}=428,\text{ }%
e_{3}^{4}=169\text{, }e_{3}^{5}=152\text{, }e_{3}^{6}=61\text{, }e_{3}^{7}=85%
\text{, }e_{3}^{8}=85\text{ and }e_{3}^{9}=85\text{.}
\end{equation*}

\textbf{Step 4}. We compute the elements $e_{4}^{i}$ to construct the matrix
$E_{i}:$%
\begin{equation*}
e_{4}^{1}=51\text{, }e_{4}^{2}=94\text{, }e_{4}^{3}=176,\text{ }e_{4}^{4}=69%
\text{, }e_{4}^{5}=63\text{, }e_{4}^{6}=25\text{, }e_{4}^{7}=35\text{, }%
e_{4}^{8}=35\text{ and }e_{4}^{9}=35\text{.}
\end{equation*}

\textbf{Step 5. }We calculate the elements $x_{i}:$%
\begin{eqnarray*}
(-1)(550) &=&51(288+5x_{1})-127(120+2x_{1}) \\
&\Rightarrow &x_{1}=2\text{.}
\end{eqnarray*}%
\begin{eqnarray*}
(-1)\left( -32\right) &=&94(192+5x_{2})-228(80+2x_{2}) \\
&\Rightarrow &x_{2}=16\text{.}
\end{eqnarray*}%
\begin{eqnarray*}
(-1)\left( 76\right) &=&176(84+5x_{3})-428(35+2x_{3}) \\
&\Rightarrow &x_{3}=5\text{.}
\end{eqnarray*}%
\begin{eqnarray*}
(-1)\left( -15\right) &=&69(24+5x_{4})-169(10+2x_{4}) \\
&\Rightarrow &x_{4}=7\text{.}
\end{eqnarray*}%
\begin{eqnarray*}
(-1)\left( -112\right) &=&63(60+5x_{5})-152(25+2x_{5}) \\
&\Rightarrow &x_{5}=12\text{.}
\end{eqnarray*}%
\begin{eqnarray*}
(-1)\left( 70\right) &=&25(204+5x_{6})-61(85+2x_{6}) \\
&\Rightarrow &x_{6}=5\text{.}
\end{eqnarray*}%
\begin{eqnarray*}
(-1)0 &=&35(60+5x_{7})-85(25+2x_{7}) \\
&\Rightarrow &x_{7}=5\text{.}
\end{eqnarray*}%
\begin{eqnarray*}
(-1)0 &=&35(60+5x_{8})-85(25+2x_{8}) \\
&\Rightarrow &x_{8}=5\text{.}
\end{eqnarray*}%
\begin{eqnarray*}
(-1)0 &=&35(60+5x_{9})-85(25+2x_{9}) \\
&\Rightarrow &x_{9}=5\text{.}
\end{eqnarray*}

\textbf{Step 6.} We rename $x_{i}$ as follows$:$%
\begin{eqnarray*}
x_{1} &=&b_{2}^{1}=2\text{, }x_{2}=b_{2}^{2}=16\text{, }x_{3}=b_{2}^{3}=5, \\
\text{ }x_{4} &=&b_{2}^{4}=7\text{, }x_{5}=b_{2}^{5}=12\text{, }%
x_{6}=b_{2}^{6}=5, \\
x_{7} &=&b_{2}^{7}=5\text{, }x_{8}=b_{2}^{8}=5\text{, }x_{9}=b_{2}^{9}=5%
\text{. }
\end{eqnarray*}

\textbf{Step 10. }We construct the block matrices $B_{i}:$%
\begin{eqnarray*}
B_{1} &=&\left[
\begin{array}{cc}
24 & 2 \\
1 & 23%
\end{array}%
\right] \text{, }B_{2}=\left[
\begin{array}{cc}
16 & 16 \\
14 & 12%
\end{array}%
\right] \text{, }B_{3}=\left[
\begin{array}{cc}
7 & 5 \\
24 & 28%
\end{array}%
\right] ,\text{ } \\
B_{4} &=&\left[
\begin{array}{cc}
2 & 7 \\
7 & 17%
\end{array}%
\right] \text{, }B_{5}=\left[
\begin{array}{cc}
5 & 12 \\
11 & 4%
\end{array}%
\right] \text{, }B_{6}=\left[
\begin{array}{cc}
17 & 5 \\
3 & 5%
\end{array}%
\right] ,\text{ } \\
B_{7} &=&\left[
\begin{array}{cc}
5 & 5 \\
5 & 5%
\end{array}%
\right] \text{, }B_{8}=\left[
\begin{array}{cc}
5 & 5 \\
5 & 5%
\end{array}%
\right] \text{, }B_{9}=\left[
\begin{array}{cc}
5 & 5 \\
5 & 5%
\end{array}%
\right] .\text{ }
\end{eqnarray*}

\textbf{Step 11.} We obtain the following message matrix $M:$%
\begin{equation*}
M=\left[
\begin{array}{cccccc}
24 & 2 & 16 & 16 & 7 & 5 \\
1 & 23 & 14 & 12 & 24 & 28 \\
2 & 7 & 5 & 12 & 17 & 5 \\
7 & 17 & 11 & 4 & 3 & 5 \\
5 & 5 & 5 & 5 & 5 & 5 \\
5 & 5 & 5 & 5 & 5 & 5%
\end{array}%
\right] =\left[
\begin{array}{cccccc}
H & A & P & P & Y & 0 \\
B & I & R & T & H & D \\
A & Y & 0 & T & O & O \\
Y & O & U & : & ) & 0 \\
0 & 0 & 0 & 0 & 0 & 0 \\
0 & 0 & 0 & 0 & 0 & 0%
\end{array}%
\right] .
\end{equation*}

\textbf{Step 9.} End of algorithm.
\end{example}

\section{Conclusion}

We have presented two new coding and decoding algorithms using the
generalized Pell $(p,i)$-numbers. Combining these algorithms together, we
can generate a new mixed algorithm such as \textquotedblleft Minesweeper
model\textquotedblright\ introduced in \cite{ucar}. Furthermore it is
possible to develop new mixed models using Fibonacci and Lucas blocking
algorithms given in \cite{Tas} and \cite{ucar}.

\end{document}